\documentclass[12pt,a4paper,reqno]{article}

\usepackage{geometry}

\usepackage[pdftitle={Lectures on Poisson Groupoids},
            pdfauthor={Camille Laurent, Mathieu Stienon, Ping Xu},
            pdfkeywords={keywords}]{hyperref}

\hypersetup{pdfpagemode=UseNone,colorlinks=false,pdfpagelayout=SinglePage,pdfstartview=FitH}

\usepackage[T1]{fontenc}
\usepackage[latin1]{inputenc}

\usepackage{setspace}
\usepackage{indentfirst}

\newcommand{\nijenhuis}{N}

\usepackage[intlimits]{amsmath}
\allowdisplaybreaks[1]
\usepackage{amssymb}
\usepackage{stmaryrd}

\usepackage{graphicx}
\usepackage[all]{xy}
\def\dar[#1]{\ar@<2pt>[#1]\ar@<-2pt>[#1]}

\usepackage{amsthm}

\theoremstyle{plain}
\newtheorem{prop}{Proposition}[section]
\newtheorem{lem}[prop]{Lemma}

\newtheorem{cor}[prop]{Corollary}
\newtheorem{thm}[prop]{Theorem}

\newtheorem*{prop*}{Proposition}
\newtheorem*{lem*}{Lemma}
\newtheorem*{sublem*}{Sublemma}
\newtheorem*{cor*}{Corollaire}
\newtheorem*{thm*}{Th\'eor\`eme}
\newtheorem*{hypo*}{Hypothesis}
\newtheorem*{conjecture*}{Conjecture}
\newtheorem*{scholum*}{Scholum}
\newtheorem{defn}[prop]{Definition}
\newtheorem*{defn*}{D\'efinition}

\theoremstyle{definition}

\newtheorem*{con*}{Construction}
\newtheorem*{note*}{Note}

\theoremstyle{remark}

\newtheorem*{shortnote*}{Note}
\newtheorem*{claim*}{Claim}
\newtheorem*{axiom*}{Axiom}

\newtheoremstyle{slanted}
  {3pt}
  {3pt}
  {\slshape}
  {}
  {\bfseries}
  {.}
  {.5em}
  {}

\theoremstyle{slanted}
\newtheorem{question}[prop]{Question}
\newtheorem*{question*}{Question}
\newtheorem{exercise}[prop]{Exercise}
\newtheorem{example}[prop]{Example}
\newtheorem{examples}[prop]{Examples}
\newtheorem*{example*}{Example}
\newtheorem*{examples*}{Examples}

\newtheorem*{ex*}{Example}
\newtheorem*{exs*}{Examples}

\newtheorem*{remark*}{Remark}
\newtheorem*{remarks*}{Remarks}
\newtheorem{rmk}[prop]{Remark}

\newtheorem*{rmk*}{Remark}
\newtheorem*{rmks*}{Remarks}
\newtheorem{warning}[prop]{Warning}
\newtheorem*{warning*}{Warning}

\DeclareMathOperator{\ad}{ad}

\newcommand{\beq}[1]{\begin{equation}\label{#1}}
\newcommand{\eeq}{\end{equation}}

\newcommand{\RR}{\mathbb{R}}

\newcommand{\derlie}[1]{\mathcal{L}_{#1}}

\newcommand{\interieur}[1]{i_{#1}} 
\newcommand{\ii}{\mathbin{\vrule width1.5ex height.4pt\vrule height1.5ex}}

\newcommand{\lie}[2]{[#1,#2]} 
\newcommand{\poisson}[2]{\{#1,#2\}} 

\newcommand{\half}{\frac{1}{2}}

\newcommand{\smalcirc}{\mbox{\tiny{$\circ$}}}

\newcommand{\gendex}[2]{\left\{ #1 \right\}_{#2}}
\newcommand{\genrel}[2]{\left\{ #1 | #2 \right\}}

\let\Vec=\overrightarrow
\let\ceV=\overleftarrow
\newcommand{\lvec}[1]{\overleftarrow{#1}}
\newcommand{\rvec}[1]{\overrightarrow{#1}}
\newcommand{\ip}[2]{\left<#1,#2\right>}
\DeclareMathOperator{\DR}{DR}
\def\lcf{\lbrack\! \lbrack}
\def\rcf{\rbrack\! \rbrack}
\newcommand{\alp }{\alpha }
\newcommand{\bet }{\beta }
\newcommand{\cala}{\mathcal{A}}
\def\gpd{\rightrightarrows}

\newcommand{\cty}{C^{\infty}}
\newcommand{\cinf}[1]{C^{\infty}(#1)}

\newcommand{\vf}{\mathfrak{X}} 
\newcommand{\XX}{\mathfrak{X}} 
\newcommand{\df}{\Omega} 

\newcommand{\toto}{\rightrightarrows}

\newcommand{\diese}{^{\sharp}}
\newcommand{\bemol}{^{\flat}}

\newcommand{\mfg}{\mathfrak{g}} 

\newcommand{\gm}{\Gamma}

\newcommand{\be }{\begin{eqnarray*}}
\newcommand{\ee }{\end{eqnarray*}}

\newcommand{\ttf}{\phi}

\newcommand{\tomega}{\tilde{\omega}}

\newcommand{\bivector}{\pi}
\newcommand{\bivectorN}{\pi_N}
\newcommand{\lieN}[2]{[#1,#2]_N}
\newcommand{\tpi}{{\tilde{\pi}}}
\newcommand{\tN}{{\tilde{N}}}

\def\til0{\Tilde{0}}

\def\dminus{\raise2pt\hbox{\vrule height1pt width 2ex}\hskip3pt}

\begin{document}

\title{Lectures on Poisson groupoids}
\author{
Camille Laurent-Gengoux \\ 
D\'epartement de math\'ematiques \\ Universit\'e de Poitiers \\ 
86962 Futuroscope-Chasseneuil, France \\ 
\href{mailto:laurent@math.univ-poitiers.fr}{\texttt{laurent@math.univ-poitiers.fr}}
\and 
Mathieu Sti\'enon
\thanks{Francqui fellow of the Belgian American Educational Foundation}
\thanks{Partially supported by the European Union through the FP6 Marie Curie R.T.N. ENIGMA
(Contract number MRTN-CT-2004-5652).} \\
E.T.H.~Z\"urich \\ Departement Mathematik \\ 
8092 Z\"urich, Switzerland \\
\href{mailto:stienon@math.ethz.ch}{\texttt{stienon@math.ethz.ch}}
\and 
Ping Xu
\thanks{Research supported by NSF grants DMS03-06665 and
DMS-0605725.}\\ 
Department of Mathematics \\ Penn State University \\
USA\\ 
\href{mailto:ping@math.psu.edu}{\texttt{ping@math.psu.edu}} }

\date{}
\maketitle

\begin{abstract}
In these lecture notes, we give a
quick account of the theory of Poisson groupoids and Lie bialgebroids.
In particular, we discuss the universal lifting theorem and its applications
including  integration of quasi-Lie bialgebroids,
integration of Poisson Nijenhuis structures and 
Alekseev and Kosmann-Schwarzbach's theory of 
$D/G$-momentum maps.
\end{abstract}

\section{Introduction}

These are the lecture notes for a mini-course given by
the third author at ICTP Trieste in July 2005.
The purpose of the mini-course was to give a
quick account of the theory of Poisson groupoids in Poisson geometry.
There are two important particular classes of Poisson groupoids: Poisson groups are one, 
the so called symplectic groupoids are another.
They correspond to two extreme cases.

Poisson groups are Lie groups which admit  compatible Poisson
structures. Poisson groups
were introduced by Drinfeld \cite{Drinfeld83a, Drinfeld83b}
as classical counterparts of quantum groups \cite{CP}
and studied by Semenov-Tian-Shansky \cite{STS}, Lu-Weinstein \cite{LuW:1990}
 and many others (see \cite{CP} for  references).
 These structures have  played an
important role in the study of integrable systems \cite{STS}.

The infinitesimal analogues of Poisson groups are
Lie bialgebras.  A Lie  bialgebra  is a pair of Lie
algebras $(\mfg,\mfg^*)$ satisfying the following compatibility
condition: the cobracket $\delta : \mfg\to
\wedge^2 \mfg$, i.e. the map dual to the Lie bracket on $\mfg^*$,
must be a Lie algebra 1-cocycle. The Jacobi identity
on $\mfg^*$ is equivalent, through dualization, to requiring that
$\delta^2=0$. Here we extend the linear operator 
 $\delta:\mfg\to\wedge^2\mfg$ to a
degree $1$ derivation $\wedge^* \mfg\to \wedge^{*+1}\mfg$
in a  natural way.  Indeed for a Lie algebra $\mfg$,
the Lie bracket on $\mathfrak{g}$
 extends naturally  to a graded Lie bracket
on $\wedge^{\bullet} \mathfrak{g} $ so that
$ (\wedge^{\bullet} \mathfrak{g}, \wedge , [\cdot,\cdot])$
is a Gerstenhaber algebra.
Using this terminology, a Lie bialgebra is nothing other than just 
a differential Gerstenhaber algebra
$(\wedge^\bullet \mathfrak{g}, \wedge , [\cdot, \cdot],\delta)$
\cite{Kosmann:1991}. A drawback of this definition is
that it does not seem obvious that the roles of $\mfg$ and
$\mfg^*$ are symmetric. A simple way to get
around this problem is to view this  differential Gerstenhaber algebra
 as one part of Kosmann-Schwarzbach's big bracket structure \cite{Kosmann:1991}.
However, a great advantage of the cobracket viewpoint 
is that the Drinfeld correspondence between  Lie 
bialgebras and Poisson groups \cite{Drinfeld83a, Drinfeld83b}
becomes much more transparent.  Namely, on the level of the integrating Lie
group $G$, the cobracket differential $\delta:  
\wedge^* \mfg\to \wedge^{*+1}\mfg$ integrates to a
multiplicative bivector field $\pi\in\mathfrak{X}^2(G)$,
and the relation $\delta^2=0$ is equivalent to
$[\pi, \pi]=0$. In other words, $(G, \pi )$ is a Poisson group.

Two decades ago, motivated by the quantization problem of
Poisson manifolds, Karasev \cite{Karasev},  Weinstein \cite{Weinstein:1987}
and later Zakrzewski \cite{Zak1, Zak2}
 introduced the notion of symplectic
groupoids. Symplectic groupoids are Lie groupoids equipped
with compatible symplectic structures. In a certain
sense, they are semi-quantum counterparts of Poisson
structures. If one thinks of Poisson manifolds as
non-linear Lie algebras, then symplectic
groupoids serve as analogues of ``non-linear Lie groups''.
And their corresponding (convolution) algebras
give rise to quantizations of
the underlying Poisson structures \cite{Weinstein:91}.
The precise link between quantization of Poisson manifolds and symplectic
groupoids was clarified only very recently by
 Cattaneo-Felder \cite{CattaneoF, CattaneoF2}.

The existence of (local) symplectic groupoids gives an
affirmative answer to a classical question: Does
every Poisson manifold admit a symplectic realization?
A symplectic realization of a Poisson manifold $(M,\pi)$ 
consists of a pair $(X,\Phi)$, where $X$ is a symplectic manifold 
and $\Phi:X\to M$ is a Poisson map which is a surjective submersion.
Hence, a symplectic realization $(X,\Phi)$  amounts  to 
 embedding  the Poisson algebra $\cty(M)$ as a Poisson subalgebra 
 of the ``symplectic algebra'' $\cty(X)$. In other words,  
a symplectic realization of a Poisson manifold $M$ desingularizes
 the Poisson structure. We refer
to \cite{Laurent} for a detailed study of this topic. 
Clearly, symplectic realizations 
are by no means unique.
The question of finding symplectic realizations of
a Poisson manifolds can be traced back to Lie \cite{Lie},
who in fact proved the local existence of such a realization
for regular Poisson structures.  The local existence for
general Poisson manifolds was obtained by Weinstein \cite{Weinstein83}
using the splitting theorem. 
And, in 1987, Karasev \cite{Karasev} and
Weinstein \cite{Weinstein:1987} independently  proved the
 existence of a global symplectic realization for any Poisson manifold.
Strikingly, they discovered that among all such realizations,
there exists a distinguished one which
 admits \emph{automatically} a local groupoid structure  compatible, in a certain sense, with the
symplectic structure. The global form of this notion is
what is now called a symplectic groupoid.
It is thus natural to explore why groupoid and symplectic structures 
arise simultaneously in such a striking
manner from a Poisson manifold.

In 1988, Weinstein introduced the notion of Poisson
groupoids, unifying both Poisson groups and
symplectic groupoids under the same umbrella \cite{Weinstein:cois}. 
It was thus tempting to develop an analogue to
Drinfeld's theory for Poisson groupoids, which
in turn could be used to study symplectic groupoids. 
In 1994, Mackenzie and one of the
authors discovered the infinitesimal counterparts of Poisson groupoids,
which are called Lie bialgebroids \cite{MackenzieX:1994}.
They can be simply characterized as differential Gerstenhaber algebras
of the form $(\gm (\wedge^\bullet A), \wedge , [\cdot, \cdot],\delta)$, where
$A$ is a vector bundle \cite{Xu:bv}. The integrability of Lie bialgebroids
was proved in \cite{MackenzieX:2000} using highly non-trivial
techniques. Given a Poisson manifold $(M,\pi) $, its cotangent bundle $T^*M$ carries a 
natural Lie algebroid structure and $(T^*M,TM)$
is indeed naturally a Lie bialgebroid. Its corresponding
differential Gerstenhaber algebra is
$(\Omega^* (M), \wedge , [\cdot, \cdot], d_{DR})$.
If $\gm\toto M$ is the $\alpha$-connected and $\alpha$-simply connected
Lie groupoid integrating the Lie algebroid structure on $T^*M$, 
then $\gm\toto M$ is a Poisson groupoid and 
the Poisson structure on $\gm$ resulting from the
integration of $d_{DR}$ turns out  to be non-degenerate. Thus one
obtains a symplectic groupoid.  We refer the reader
to \cite{CF, CF2} for integrability criteria for Poisson manifolds. 
 Another proof of the existence of symplectic   groupoids
was recently obtained by Cattaneo-Felder \cite{CattaneoF}
using the Poisson sigma model. 

Quasi-Lie bialgebroids are the infinitesimal objects associated to 
quasi-Poisson groupoids. They were introduced by Roytenberg \cite{Roy}
as a generalization of Drinfeld's quasi-Lie bialgebras \cite{Drinfeld}.
Twisted Poisson structures \cite{SW}
 are examples of quasi-Lie bialgebroids \cite{Roy}. 
It is thus a natural question to ask  whether every quasi-Lie bialgebroid
integrates to a quasi-Poisson groupoid. 
On the other hand, quasi-Poisson groupoids also arise naturally in the study
of generalized momentum map theory. In \cite{AK-S},    
Alekseev and Kosmann-Schwarzbach introduced \emph{quasi-Poisson spaces}
with $D/G$-momentum maps, which are generalizations of
quasi-Hamiltonian spaces with group valued momentum maps
\cite{AMM} (see also \cite{LX, Xu:JDG04} from the view point of
Hamiltonian $\gm$-spaces). It turns out that these quasi-Poisson spaces
are exactly Hamiltonian $\gm$-spaces, where $\gm$ denotes a 
quasi-Poisson groupoid. A primordial tool for integrating 
quasi-Lie bialgebroids is the so called 
\emph{universal lifting theorem}: for an $\alpha$-connected and
$\alpha$-simply connected Lie groupoid $\gm$ there is a natural
isomorphism between the graded Lie algebra of multiplicative
multi-vector fields on $\gm$ and the the graded Lie algebra
of multi-differentials on its Lie algebroid $A$ --- see Section~\ref{4.1} 
for the precise definition of multi-differentials.
Many well-known integration theorems (in particular 
of Lie bialgebroids and of twisted Poisson manifolds) 
are easy consequences of this universal lifting theorem. 
This is also the viewpoint we take in these notes.
In particular, we discuss the integration problem of
Poisson Nijenhuis structures as an application.

The notes are organized as follows.  In Section 2, we  give a
quick account of Poisson group and Lie bialgebra theory.
In Section 3, we introduce Poisson groupoids and Lie bialgebroids.
Section 4 is devoted to the universal lifting theorem and
its applications including quasi-Poisson groupoids
and integration of Poisson quasi-Nijenhuis structures.

{\bf Acknowledgments.}
We would like to thank the organizers of the ICTP program
on ``Poisson geometry'' for organizing such a wonderful event.


\section{Poisson groups and Lie bialgebras}

\subsection{From Poisson groups to Lie bialgebras}

Let us recall the definition of Poisson groups.

\begin{defn}[\cite{Drinfeld83a, Drinfeld83b}]
A Poisson group is a Lie group endowed with
a Poisson structure $\pi \in {\mathfrak X}^2 (G)$ such that the  multiplication
$m: G \times G \to G$ is a Poisson map, where $G \times G$
is equipped with the product Poisson structure. 
\end{defn}

\begin{example}
The reader may have in mind the following  two trivial examples.
\begin{enumerate}
\item For any Lie algebra ${\mathfrak g}$,
its dual (${\mathfrak g}^* ,+$) is a Poisson group
where {\it (i)} the Lie group structure is given by 
the addition
{\it (ii)} the Poisson structure is the linear Poisson structure, i.e. Lie-Poisson structure.
\item Any Lie group $G$ is a Poisson group
with respect to the trivial Poisson bracket.
\end{enumerate}
\end{example}

To impose that $m : G \times G \to G$
is a Poisson map is equivalent to
impose any of the following two  conditions
\begin{enumerate}
\item $\forall g,h \in G$, $  m_* (\pi_g , \pi_h) = \pi_{gh} \, \,$ or,
\item $\forall g,h \in G$, $(R_h)_* \pi_g + (L_g)_* \pi_h = \pi_{gh} $.
\end{enumerate}
(Here $\pi_p \in \wedge^2 T_p G $ stands for  the value of the bivector
field at the point $p \in P$.)

This leads to the following definition:

\begin{defn}[\cite{Lu:1990}]
A bivector field $\pi $ on $G$ is said to be multiplicative
if 
\begin{equation}
\label{eq:mult}
 (R_h)_* \pi_g + (L_g)_* \pi_h = \pi_{gh} , \ \ \ \forall g, h \in G.
\end{equation}
\end{defn}

In particular, $\pi \in \mathfrak{X}^2 (G)$ endows $G$ with a
 structure of Poisson group
if, and only if, {\it (i)}  the identity $[\pi,\pi]=0$ holds
and {\it (ii)} $\pi $ is multiplicative.

\begin{rmk}
Any multiplicative bivector $\pi$
vanishes in $g=1$, where $1$ is the unit element of the group $G$.
This can be seen from Eq.~\eqref{eq:mult} by letting $g=h=1$.
\end{rmk}

It is sometimes convenient to consider 
$\tilde{\pi}(g) = (R_g)^{-1}_* \pi_g $
rather that $\pi $ itself. Note that $\tilde{\pi}(g) $
  is, by construction, a smooth map from
$G$ to $\wedge^2 \mathfrak{g}$
(where, implicitly, we have identified the Lie algebra
 $\mathfrak{g}$ with the tangent space at $g=1$ of the Lie group $G$).
When written  with the help  of $\tilde{\pi}$,
the condition $(R_h)_* \pi_g + (L_g)_* \pi_h = \pi_{gh}$
reads 
\begin{align*} (R_{gh})^{-1}_*  [(R_h)_* \pi_g + (L_g)_* \pi_h ] &= (R_{gh})^{-1}_* \pi_{gh} \\   \tilde{\pi}(g) + Ad_g \tilde{\pi}(h) &= \tilde{\pi}(gh) .\end{align*}
I.e. $ \tilde{\pi}: G \to \wedge^2 \mathfrak{g} $
is a Lie group $1$-cocycle, where $G$ acts on $ \wedge^2 \mathfrak{g}$ by
adjoint action.

Now, differentiating a Lie group $1$-cocycle
at the identity, one gets a Lie algebra
$1$-cocycle $  \mathfrak{g} \to \wedge^2 \mathfrak{g}$. 
For example, the  $1$-cocycle 
$\delta: \mathfrak{g} \to \wedge^2 \mathfrak{g} $ associated to 
the Poisson structure $\tilde{\pi}$ is given by 
\begin{multline*} 
\delta(X) =\frac{d}{dt}_{|_{t=0}}\tilde{\pi}\big(\exp(tX)\big) 
=\frac{d}{dt}_{|_{t=0}}(R_{\exp(-tX)})_*\pi_{\exp(tX)} \\
=(\phi_{-t})_* \pi_{\phi_t (1) }
=(L_{\ceV{X}} \pi)_{|_{g=1}} ,
\end{multline*}
where $X$ denotes any element of $\mathfrak{g}$, $\ceV{X}$ is the left invariant vector field
on $G$ corresponding to $X$ and $\phi_t$ is its flow.

We have therefore determined $L_{\ceV{X}} \pi$ at $g=1$.
We now try to compute it 
at other points.

For all $g \in G$, since $\pi$ is multiplicative,
we have $\forall X\in \mathfrak{g}$,
\begin{align*} \pi_{ g \, \exp(tX) }&=  (R_{\exp(-tX)})_* \pi_g + (L_g)_* \pi_{\exp(tX)} \\   
  (R_{\exp(-tX)})_*  \pi_{g \, \exp(tX)}  &=   
 \pi_g + (R_{\exp(-tX)})_* (L_g)_* \pi_{\exp(tX)} \\ 
  (\phi_{-t})_* \pi_{ \phi_t (g)} &=  \pi_g + L_g (\phi_{-t})_* \pi_{\phi_t (1)} .\end{align*}
Taking the derivative of the previous identity at  $t=0$,
one obtains:
\[ (L_{\ceV{X}} \pi)_{|_g} = (L_g)_*  L_{\ceV{X}}\pi_{|1} = (L_g)_* \delta (X) .\]
which implies that $L_{\ceV{X}} \pi$ is left invariant.
For all $Y \in \wedge^k \mathfrak{g}$, we denote 
by $ \ceV{Y} $ (resp. $\Vec{Y}$) the left (resp. right)
invariant $k$-vector field  on $G$ equal to $Y$
at $g=1$.  Then we obtain 
the following formula:
\[ L_{\ceV{X}} \pi = \ceV{\delta (X)} .\]
and, for similar reasons
\[ L_{\Vec{X}} \pi = \Vec{\delta (X)} .\]
We can now extend $\delta:\mathfrak{g} \to\wedge^2 \mathfrak{g} $
to a derivation   of degree $+1$ of the graded commutative associative
algebra $\wedge^* \mathfrak{g} $ that we denote by the same symbol
 $\delta:\wedge^{\bullet} \mathfrak{g} \to\wedge^{\bullet+1} \mathfrak{g} $.

\begin{lem} 
\begin{enumerate}
\item  $\delta^2=0$
\item $\delta [X, Y]=[\delta X, Y]+[X, \delta  Y], \ \ \forall X, Y\in \mathfrak{g} $
\end{enumerate}
\end{lem}

\begin{proof}
(1) For all $X \in \mathfrak{g}$,
\[ [ \ceV{X}, [\pi,\pi]] = 2[[\ceV{X},\pi],\pi]
= 2[\ceV{\delta (X)},\pi] = 2\ceV{\delta^2(X)} .\]
But $[\pi,\pi]=0$, hence $ \delta^2 (X)=0$.

(2) This follows from the graded Jacobi identity:
\[ [\ceV{[X, Y]}, \pi]=[[\ceV{X}, \ceV{Y}], \pi]
=[[\ceV{X}, \pi ], \ceV{Y}]+[\ceV{X}, [\pi , \ceV{Y}]] . \qedhere \]
\end{proof}

The infinitesimal object associated to a Poisson-Lie group is
therefore as defined below. 

\begin{defn}
A Lie bialgebra is a Lie algebra $ \mathfrak{g}$
equipped with a  degree $1$-derivation $\delta$
of the graded commutative associative algebra $ \wedge^\bullet \mathfrak{g}$
 such that 
\begin{enumerate}
\item $\delta([X,Y])=  [\delta(X),Y]+ [X, \delta(Y)]$ and
\item $\delta^2 =0$.
\end{enumerate}
\end{defn}

\begin{rmk}
Recall that a \emph{Gerstenhaber algebra} $A = \oplus_{i \in {\mathbb N}} A^i$
is a graded commutative algebra s.t. $A = \oplus_{i \in {\mathbb N}} A^{(i)} $
(where $A^{(i)}=A^{i+1}$) is a graded Lie algebra with the compatibility condition
\[ [a,bc]=[a,b]c +  (-1)^{(|a|+1) |b|} b [a,c] \] 
for any $a \in A^{|a|}$, $b\in A^{|b|}$ and $c\in A^{|c|}$.

A \emph{differential Gerstenhaber algebra} is a Gerstenhaber algebra equipped with a degree 1
derivation of square zero and compatible with respect to both brackets \cite{Xu:bv}.

The Lie bracket on
$\mathfrak{g}$
can be extended to a graded Lie bracket
on $\wedge^{\bullet} \mathfrak{g} $ so that  $ (\wedge^{\bullet} \mathfrak{g}, \wedge , [\cdot,\cdot])$
is a Gerstenhaber algebra.
Using this terminology, a 
Lie bialgebra is nothing else than a differential Gerstenhaber algebra 
$(\wedge^\bullet \mathfrak{g}, \wedge , [\cdot,\cdot],\delta)$.
\end{rmk}

Given a Lie bialgebra $(\mathfrak{g},\delta)$, let us consider its dual
 $\delta^*: \wedge^2 \mathfrak{g}^* 
\to \mathfrak{g}^*$ of the derivation $\delta$.

Let $ [\xi,\eta]_{\mathfrak{g}^*}= \delta^* (\xi \wedge \eta)$
for all $\xi,\eta \in \mathfrak{g}^*$.
The bilinear map $(\xi,\eta) \to [\xi,\eta]_{\mathfrak{g}^*} $ is 
skew-symmetric and 
\[  \delta^2=0   \hspace{.5cm}  \Leftrightarrow  \hspace{.5cm}  [\cdot,\cdot]_{\mathfrak{g}^*}  \text{ satisfies the Jacobi identity} \]

Therefore, the dual $\mathfrak{g}^*$ of a Lie bialgebra $(\mathfrak{g}, \delta)$ is a Lie algebra again  (which justifies the name).
Conversely, a Lie bialgebras can be described again by: 

\begin{prop} (see \cite{CP} Chapter 1 for instance)
A Lie  bialgebra is  equivalent to a pair of Lie algebras
$(\mathfrak{g}, \mathfrak{g}^*)$ compatible in the following sense:
 the coadjoint action  of $\mathfrak{g}$ on $\mathfrak{g}^*$
is a derivation of the bracket  $ [\cdot,\cdot]_{\mathfrak{g}^*} $,
i.e.
\[  \ad^*_X [\alpha,\beta]_{\mathfrak{g}^*}= [\ad^*_X \alpha,\beta]_{\mathfrak{g}^*}
+   [\alpha,\ad^*_X \beta]_{\mathfrak{g}^*}\qquad \text{for all $X \in \mathfrak{g}, \alpha, \beta \in \mathfrak{g}^*$} .\]
\end{prop}

\begin{rmk} Note that Lie bialgebras are in duality: namely
$(\mathfrak{g}, \mathfrak{g}^*)$ is a Lie bialgebra if and only
if $(\mathfrak{g}^*, \mathfrak{g})$ is a Lie bialgebra.
This picture can be seen more naturally using
Manin triples, which will be discussed in the next lecture.
\end{rmk}

\subsection{$r$-matrices}

We now turn our  attention to a particular class of Lie bialgebras, i.e.
those coming from $r$-matrices.

We start from a Lie algebra $\mathfrak{g}$.
Assume that we are given an element $r\in \wedge^2 \mathfrak{g} $.
Then define $\delta$ by, for all $X \in \wedge^\bullet \mathfrak{g}$,
$\delta (X) = [r,X]$. As can easily be checked, $\delta$ is a derivation  of 
$\wedge^\bullet \mathfrak{g}$. Note that, in terms of (Chevalley-Eilenberg) 
cohomology, $ \delta$ is the coboundary of $r \in \wedge^2 \mathfrak{g}$.

The condition $\delta^2 (X) =0$ is equivalent to
the relation $[X,[r,r]]=0$,
which itself holds if, and only if, $[r,r]$ is $\ad$-invariant.
Conversely, any $r\in \wedge^2 \mathfrak{g} $
such that  $ [r,r] $ is  $\ad$-invariant defines a Lie bialgebra.
Such an $r$ is called an {\it $r$-matrix}.
If moreover $[r,r]=0$, then this Lie bialgebra is called \emph{triangular}.

Here are two well-known examples of $r$-matrices.
\begin{example} [\cite{Drinfeld83a, Drinfeld83b, STS, CP}]
\label{1.2} 
\begin{enumerate}
\item Consider a semi-simple Lie algebra $\mathfrak{g}$ of rank $k$ over $\mathbb{C}$ and a Cartan sub-algebra $\mathfrak{h}$.
Let $\{ e_\alpha, f_\alpha,  \alpha \in \Delta_+\} \cup \{ h_i , i =1,\dots, k  \}$
be a Chevalley basis.
Then $ r  = \sum_{\alpha \in \Delta_+} \lambda_{\alpha} e_\alpha \wedge f_\alpha $,
where $\lambda_{\alpha} = \frac{1}{(e_\alpha,f_\alpha)}$,
is an $r$-matrix.
\item Consider now a compact semi-simple Lie algebra
$\mathfrak{k}$ over $\mathbb{R}$.
Let $\{ e_\alpha, f_\alpha,  \alpha \in \Delta_+ \} \cup \{ h_i, i =1\dots, k \}$
be a Chevalley basis (over $\mathbb{C}$) of 
the complexified Lie algebra $\mathfrak{g}= \mathfrak{k}^{\mathbb{C}}$, that we assume to be constructed
so that the family $\{ X_\alpha, Y_\alpha, \alpha \in \Delta_+\} \cup \{ 
t_i , i =1,\dots,k  \} $ is a basis of $\mathfrak{k}$ (over $\mathbb{R}$) where
\[
 \left\{ \begin{aligned} X_\alpha &= e_\alpha - f_\alpha && \text{for all $\alpha \in {\Delta}_+ $} \\     
 Y_\alpha &= \sqrt{-1} (e_\alpha + f_\alpha) &&
 \text{for all $\alpha \in {\Delta}_+ $} \\ 
 t_i &= \sqrt{-1} h_i &&
 \text{for all $ i \in \{1,\dots,k\}$} . 
  \end{aligned} \right. \]
Let \[ \hat{r} = \sqrt{-1} \, r  = \sqrt{-1} \, \sum_{\alpha \in \Delta_+} \lambda_{\alpha} e_\alpha \wedge f_\alpha .\] Then
$\hat{r}$ is, according to the first example above,
 an $r$-matrix of $\mathfrak{g}=\mathfrak{k}^{\mathbb{C}}$. However, by a direct computation, one checks that
\[ \hat{r} =  \frac{1}{2} \sum_{\alpha \in \Delta_+} \lambda_{\alpha} X_\alpha \wedge Y_\alpha \] 
so that $\hat{r} $ is indeed an element of
$ \wedge^2 \mathfrak{k}$, and therefore is
 an $r$-matrix on $\mathfrak{k}$.
Hence, it defines a Lie bialgebra
structure on the real Lie algebra $\mathfrak{k}$.
\end{enumerate}
\end{example}

\subsection{Lie bialgebras and simply-connected Lie groups}

We have already explained how to
get a Lie bialgebra from a Poisson group.
The inverse  is true as well
when the Lie group is connected and simply-connected.

\begin{thm}[Drinfeld \cite{Drinfeld83a}] 
Assume that $G$ is a connected
and simply-connected Lie group. Then there exists
a one-to-one correspondence
   between Poisson groups $(G,\pi) $ and Lie bialgebras $ (\mathfrak{g},\delta)$.
\end{thm}

\begin{example}
In particular, for a Lie bialgebra coming from
an $r$-matrix $r$, the corresponding Poisson structure on $G$
is the bivector field  $\ceV{r}- \Vec{r} $.
\end{example}

Applying the theorem above to the previous two examples, we are
lead to

\begin{prop}
\begin{enumerate}
\item \cite{Drinfeld83a, Drinfeld83b, CP} Any complex semi-simple Lie group admits
a natural (complex)  Poisson group structure. 
\item  \cite{Soibel, LuW:1990} Any compact semi-simple Lie group admits
a natural Poisson group structure,
called the Bruhat-Poisson structure.
\end{enumerate}
\end{prop}

\begin{rmk} Poisson groups come in pairs in the following sense. Given a Poisson group $(G,\pi)$, let $(\mathfrak{g},\mathfrak{g}^*)$ be its Lie bialgebra, then we know $(\mathfrak{g}^*,\mathfrak{g})$ is also a Lie bialgebra which gives rise to a Poisson group denoted $(G^*,\pi')$. \end{rmk}

\begin{example}
[\cite{LuW:1990, Lu:1990}]
Since any element $g$ of the Lie group $G=SU(2)$ is of the form,
\[ g = \left( \begin{array}{cc} \alpha & \beta \\ 
-\bar{\beta} & \bar{\alpha} \\ \end{array} \right) ,\]
we can define complex coordinate functions $\alpha$ and $\beta$ on $G$.
Note that these coordinates are not ``free'' since $|\alpha|^2+|\beta|^2 =1$.
The Bruhat-Poisson structure is given by
\begin{align*}  \{\alpha, \bar{\alpha} \}  &= 2 \sqrt{-1} \, \beta \bar{\beta} & 
\{\alpha, \beta \} &= - \sqrt{-1} \, \alpha \beta  \\ 
\{\alpha, \bar{\beta}\} &= - \sqrt{-1} \, \alpha \bar{\beta} &
\{\beta, \bar{\beta}\} &= 0 
.\end{align*} 
\end{example}

\newcommand{\genrelamoi}[2]{\left\{ #1 \right| \left. #2 \right\}}

\begin{example}
Below are two examples of duals of Poisson groups.
\begin{enumerate}
\item For the Poisson group $G= SU(2)$ equipped
with the Bruhat-Poisson structure, the dual group
is  \cite{LuW:1990, Lu:1990}
\[ G^* = SB(2) \simeq \gendex{ \left. \left( \begin{array}{cc}  a  & b+\sqrt{-1} c \\    0 & \frac{1}{a} \end{array}  \right) \right| \begin{array}{c} \\ \, b,c \in {\mathbb R}, a \in {\mathbb R}^+  \\  \\ \end{array} }{} 
\]
Using these coordinates, the Poisson structure on $G^*$ is given explicitly by
\begin{align*}   \{b,c\}& =  a^2 - \frac{1}{a^2}   \\ \{a,b\}   & = ab  \\ \{a,c\} & =  ac   .\end{align*}
\item For the Poisson group $G= SL_{{\mathbb C}}(n)$,  equipped
with the Poisson bracket constructed in Example~\ref{1.2}, the dual group
is \cite{STS} 
\[ \begin{array}{ccccc}G^*  &= &B_+ \star B_-    & \simeq &  \genrelamoi{\begin{array}{c} \\ (A,B)\\   \\ \end{array}}{  
 \begin{array}{c}  A \mbox{ upper triangular with determinant $1$, } \\ B \mbox{ lower triangular with determinant $1$},  \\ \mbox{s.t.} \, \, \,   \mbox{diag}(A) \cdot \mbox{diag}(B)=1 \end{array}}  \end{array}
\]
\end{enumerate}
\end{example}

\subsection{Poisson group actions}

\begin{defn}[\cite{STS}]
Let $G$ be a Poisson group.
Assume that $G$ acts on a Poisson manifold $X$.
The action is said to be a Poisson action
if the action map
\[ G \times X \to X:(g,x)\mapsto g \cdot x \] 
is a Poisson map, where $G \times X$ is equipped
with the product Poisson structure.
\end{defn}

\begin{warning} Note that in general, this definition \emph{does not}
imply that,
for a fixed $g$ in $G$, the action $x \mapsto g \cdot x$
is a Poisson automorphism of $X$. 
The reader should not confuse Poisson actions
with actions preserving the Poisson structure!
Note, however, that when the Poisson
structure on the Lie group $G$ is the trivial one,
then a Poisson action is an action of $G$ on $X$ which preserves the Poisson structure.
\end{warning}

\begin{example}
Any Lie group $G$ acts on itself by  left translations.
If $G$ is a Poisson group then this action is a Poisson action.
\end{example}

\begin{prop}[\cite{Lu:1990}]
Let $G$ be a Poisson group with Lie bialgebra
 $(\mathfrak{g},\delta) $.
Assume that $G$ acts on  a  manifold $X$
and let $ \rho: \mathfrak{g} \to \mathfrak{X}^1(X) $
be the infinitesimal action.
The action of $G$ on $X$ is a Poisson action
if and only if the following diagram commutes
\[ \xymatrix{ \mathfrak{g} \ar[r]^{\rho} \ar[d]_{\delta} & 
\mathfrak{X}(X) \ar[d]^{[\pi ,\cdot]} \\ 
\wedge^2\mathfrak{g} \ar[r]_{\rho} & \mathfrak{X}^2(X) } \]  
\end{prop}

In terms of Gerstenhaber algebras,
the commutativity of the previous diagram has a clear meaning:
it simply means that $\rho: \wedge^{\bullet} \mathfrak{g} \to \mathfrak{X}^{\bullet}(X)$ 
is a morphism of differential Gerstenhaber algebra.
In other words, the natural map $T^*X\to  \mathfrak{g}^*$ induced by
the infinitesimal $ \mathfrak{g} $-action is  a Lie bialgebroid morphism \cite{Xu:1995}.

\begin{example}
For the dual $SL_{\mathbb{C}}(3)^* =  B_+ \star B_- $ of $G=SL_{\mathbb{C}}(3)$.
Consider the Poisson manifold
\[  X=  \genrelamoi{  \left( \begin{array}{ccc}  1  &  x & y  \\ 0 & 1 & z \\
    0 & 0 & 1 \end{array}  \right) }{  \begin{array}{c} \\ x,y,z \in \mathbb{C} \\  \\ \end{array} } \]
equipped with the Poisson bracket \cite{Dubrovin}
\begin{align*} 
\{x,y\} &= xy-2z \\ 
\{y,z\} &= yz-2x \\ 
\{z,x\} &= zx-2y \\ 
\end{align*}

The Lie group $ G^*= B_+ \star B_-$ acts
on $X$ by 
\[ (A,B) \cdot U \mapsto A U B^{T} \]
with $(A,B) \in B_+ \star B_- \simeq G^*$ and $U \in X$.
This action turns to be a Poisson action \cite{Xu:Dirac}.
\end{example}


\section{Poisson groupoids and Lie bialgebroids}

\subsection{Poisson and symplectic  groupoids}

In this section,
we first introduce the notion of Poisson groupoids.
The notion of Poisson groupoids was introduced by
Weinstein \cite{Weinstein:cois} as
a unification of  Poisson groups and symplectic
groupoids:
\[ \text{Poisson groupoids } \begin{cases} \text{ Poisson groups} \\ \text{ symplectic groupoids} \end{cases} \]

Let $\Gamma\rightrightarrows M$ be a Lie groupoid. A Poisson groupoid structure on $\Gamma$ should be a multiplicative Poisson structure on $\Gamma$.

To make this more precise, recall that
 in the  group case, the following are equivalent
\begin{align*} & \pi \text{ is multiplicative} \\
\Leftrightarrow \quad  & m:G\times G\to G \text{ is a Poisson map} \\
\Leftrightarrow \quad & \genrel{(x,y,xy)}{x,y\in G}\subset G\times G\times \overline{G} \text{ is coisotropic}
\end{align*}
where $\overline{G}$ denotes $(G,-\pi)$.
This motivates the following definition.

\begin{defn}
\cite{Weinstein:cois}
A groupoid $\Gamma$ with a Poisson structure $\pi$ is said to be a Poisson groupoid if the graph of the groupoid multiplication
\[ \Lambda=\genrel{(x,y,xy)}{(x,y)\in \Gamma_2 \text{ composable pair}}\subset \Gamma\times\Gamma\times\overline{\Gamma} \] 
is coisotropic. Here $\overline{\Gamma}$ means that $\Gamma$ is equipped with the opposite Poisson structure $-\pi$.
\end{defn}

\begin{examples}
\begin{enumerate}
\item If $P$ is a Poisson manifold, then $P\times\overline{P}\rightrightarrows P$ is a Poisson groupoid \cite{Weinstein:cois}.
\item Let $A$ be the Lie algebroid of a Lie groupoid $\Gamma$ and $\Lambda\in\Gamma(\wedge^2 A)$ be an element satisfying $\mathcal{L}_X[\Lambda,\Lambda]=0$, $\forall X\in\Gamma(A)$. Then $\pi=\lvec{\Lambda}-\rvec{\Lambda}$ defines a Poisson groupoid structure on $\Gamma$ \cite{LiuXu:1996}.
\end{enumerate}
\end{examples}

\begin{defn} [\cite{CDW}]
A symplectic groupoid is a Poisson groupoid $(P\rightrightarrows M,\pi)$ such that $\pi$ is non-degenerate. In other words, $\Lambda\subset\Gamma\times\Gamma\times\overline{\Gamma}$ is a Lagrangian submanifold.
\end{defn}

\begin{examples}
\begin{enumerate}
\item $T^*M\rightrightarrows M$ with the canonical cotangent symplectic structure is a symplectic groupoid
\item If $G$ is a Lie group, then $T^*G\rightrightarrows
  \mathfrak{g}^*$ is a symplectic groupoid \cite{CDW}. 
Here the symplectic structure on $T^*G$ is the canonical cotangent symplectic structure. The groupoid structure is as follows. Right translations give an isomorphism between $T^*G$ and the transformation groupoid $G\times \mathfrak{g}^*$ where $G$ acts on $\mathfrak{g}^*$ by coadjoint action.
\item In general, if $\Gamma\rightrightarrows M$ is a Lie groupoid with Lie algebroid $A$, then $T^*\Gamma\rightrightarrows A^*$ is a symplectic groupoid.
Here the groupoid structure can be described as follows.
Let $\Lambda\subset \Gamma\times\Gamma\times\Gamma$ denote the graph
 of the multiplication and 
$N^*\Lambda\subset T^*\Gamma\times T^*\Gamma\times T^*\Gamma$ 
its conormal space.
One shows that $\overline{N^*\Lambda}=\genrel{(\xi,\eta,\delta)}{(\xi,\eta,-\delta)\in N^*\Lambda}$
 is the graph of a groupoid multiplication on $T^*\Gamma$
 with corresponding unit space isomorphic to $A^*\simeq N^*M$.
This defines a groupoid structure on $T^*\Gamma\rightrightarrows A^*$.
\end{enumerate}
\end{examples}

Symplectic groupoids were introduced by
  Karasev \cite{Karasev},   Weinstein \cite{Weinstein:1987}, 
and  Zakrzewski \cite{Zak1, Zak2} in their  study of deformation
quantization of Poisson manifolds. Their relevance
with the star products was  shown by the work of
Cattaneo-Felder \cite{CattaneoF, CattaneoF2}.
On the other hand, the existence of  (local)
symplectic groupoids solves a classical
question regarding  symplectic realizations of Poisson manifolds \cite{CDW},
which can be described as follows.
Given a Poisson manifold $M$,  is it possible to
embed the Poisson algebra $\cty(M)$ into a Poisson subalgebra of 
$\cty(X)$,  where $ X$ is a symplectic manifold?
Note that  according to Darboux theorem
 there exist local coordinates $(p_1,\dots,p_k,q_1,\dots,q_k)$ in
 which the Poisson bracket on $\cty(X)$ has the following form:
 $\poisson{p_i}{q_j}=\delta_{ij}$, $\poisson{p_i}{p_j}= \poisson{q_i}{q_j}=0$.
Hence locally this amounts to finding independent functions
$\Phi_i (p_1,\dots,p_k,q_1,\dots,q_k), \ i=1, \cdots , r$ such that
$\{\Phi_i, \Phi_j\}=\pi_{ij} (\Phi_1, \cdots, \Phi_r)$, where
the left hand side stands for the Poisson bracket in $\RR^{2k}$ with 
respect to the Darboux coordinates and $\sum_{ij}\pi_{ij}\frac{\partial}{
\partial x_i}\frac{\partial}{ \partial x_j}$
 is the  Poisson tensor on $M$.
This is exactly what  Lie  first investigated in 1890  under
 the name of ``Function groups'' \cite{Lie}.
Let us give a  precise definition below.

\begin{defn}
A symplectic realization of a Poisson manifold $(M,\pi)$ consists of a pair $(X,\Phi)$, where $X$ is a symplectic manifold and $\Phi:X\to M$ is a Poisson map which is a surjective submersion.
\end{defn}

This leads to the following natural

\begin{question} Given a Poisson manifold, does 
there exist a symplectic realization? And if so, is it unique?
\end{question}

The  local existence of symplectic realizations was
proved by Lie in the constant rank case \cite{Lie}.
For a general Poisson manifold, it was proved
by Weinstein \cite{Weinstein83} in 1983  using the splitting theorem.

It is clear that  symplectic realizations of a given
Poisson manifold are not unique. The following
Karasev-Weinstein theorem states that  there
exist a canonical global  symplectic realization for
any Poisson manifold.

\begin{thm}[\cite{Karasev, Weinstein:1987}]
\begin{enumerate}
\item Symplectic realizations exist globally for any Poisson manifold;
\item Among all symplectic realizations, there exists
 a distinguished symplectic realization, which admits a
 compatible local groupoid structure making it into  a
  symplectic local groupoid.
\end{enumerate}
\end{thm}

The original proofs are highly non-trivial. The  idea  was
to use local symplectic realizations and    to patch them together.
A  different proof was recently obtained by  Cattaneo-Felder using Poisson
sigma models \cite{CattaneoF}.

Another approach due to Mackenzie and one
of the authors is to consider  integrations of
Lie bialgebroids \cite{MackenzieX:1994, MackenzieX:2000}.
An advantage of this approach is the clarification of why symplectic and groupoid structures arise in 
the context of Poisson manifolds in such a striking manner.
 
The following theorem gives an equivalent characterization
of Poisson groupoids.

\begin{thm}[\cite{Xu:1995}]
\label{thm:Poissongp}
Let $\Gamma\overset{\alpha}{\underset{\beta}\rightrightarrows}M$ be a Lie groupoid. Let $\pi\in\mathfrak{X}^2(\Gamma)$ be a Poisson tensor. Then $(\Gamma,\pi)$ is a Poisson groupoid if and only if all the following hold.
\begin{enumerate}
\item For all $(x,y)\in\Gamma_2$, \[ \pi(xy)=R_Y\pi(x)+L_X\pi(y)-R_Y L_X\pi(w) ,\] where $w=\beta(x)=\alpha(y)$ and $X, Y$ are (local) bisections through $x$ and $y$ respectively.
\item $M$ is a coisotropic submanifold of $\Gamma$
\item For all $x\in\Gamma$, $\alpha_*\pi(x)$ and $\beta_*\pi(x)$ only depend on the base points $\alpha(x)$ and $\beta(x)$ respectively.
\item For all $\alpha,\beta\in\cty(M)$, one has $\poisson{\alpha^*f}{\beta^*g}=0$, $\forall f,g \in\cty(M)$.
\item The vector field $X_{\beta^*f}$ is left invariant for all $f\in\cty(M)$.
\end{enumerate}
\end{thm}

\begin{rmk}
If $M$ is a point, then (2) implies that $\pi(1)=0$, which
together with (1) implies the multiplicativity condition.
It is easy to see that (3)-(5) are  automatically satisfied.
Thus one obtains the characterization of a Poisson group: a Lie group equipped with a multiplicative Poisson tensor.
\end{rmk}

\subsection{Lie bialgebroids}

In order to study  the infinitesimal counterparts
of Poisson groupoids, we follow the situation of Poisson groups.
As a consequence of Theorem~\ref{thm:Poissongp},
we have  the following:

\begin{cor} 
\label{corollary1}
Given a Poisson groupoid $(\Gamma\rightrightarrows M,\pi)$, we have
\begin{enumerate}
\item for all $X\in\Gamma(A)$, $[\lvec{X},\pi]$ is still left invariant
\item $\pi_M := \alpha_*\pi$ (or $-\beta_*\pi$) is a Poisson tensor on $M$
\end{enumerate}
\end{cor}
\begin{proof}
For all $X\in\Gamma(A)$, take $\xi_t=\exp tX\in U(\Gamma)$ (the space of bisections of $\Gamma$), $u_t=(\exp tX)(u)$ and $x\in \Gamma$ with $\beta(x)=u$.
In other words, $u_t$ is the flow of $\lvec{X}$ initiated at $u$.
Let $K$ be any bisection through $x$.
One gets
\begin{align*} & \pi(xu_t)=R_{\xi_t}\pi(x)+L_K\pi(u_t)-L_K R_{\xi_t} \pi(u) \\
\Rightarrow \quad & R_{{\xi_t}^{-1}}\pi(xu_t)=\pi(x)+L_K R_{{\xi_t}^{-1}} \pi(u_t)- L_K\pi(u) \in \wedge^2 T_x\Gamma
\end{align*}
and, differentiating with respect to $t$ at $0$,
\[ (\mathcal{L}_{\lvec{X}}\pi)(x)=L_K\big((\mathcal{L}_{\lvec{X}}\pi)(u)\big) .\]
This implies that $\mathcal{L}_{\lvec{X}}\pi$ is left invariant.
\end{proof}

Now, we can introduce operators
 $\delta:\Gamma(\wedge^i A)\to\Gamma(\wedge^{i+1} A)$.
For $i=0$, \[ \cty(M)\to\Gamma(A):f\mapsto X_{\beta_*f}=[\beta^* f,\pi] .\]
For $i=1$, \[ \Gamma{A}\to\Gamma(\wedge^2 A):X\mapsto \lvec{\delta X}=[\lvec{X},\pi] .\]
The following lemma can be easily verified.

\begin{lem}
\begin{enumerate}
\item \label{f1} $\delta(fg)=g\delta f+f\delta g, \quad \forall f,g\in\cty(M)$
\item \label{f2} $\delta(fX)=\delta f\wedge X+f\delta X, \quad \forall f\in\cty(M) \text{ and } X\in\Gamma(A)$
\item \label{f3} $\delta\lie{X}{Y}=\lie{\delta X}{Y}+\lie{X}{\delta Y}, \quad \forall X,Y\in\Gamma(A)$
\item \label{f4} $\delta^2=0$
\end{enumerate}
\end{lem}

\begin{defn}
A Lie bialgebroid is a Lie algebroid $A$ equipped with a degree 1 derivation $\delta$ of the associative algebra $(\Gamma(\wedge^{\bullet}A),\wedge)$ satisfying conditions \ref{f3} and \ref{f4} of the previous lemma.
\end{defn}

\begin{exercise}
Show that a Lie bialgebroid structure is equivalently characterized as a degree 1 derivation $\delta$ of the Gerstenhaber algebra $(\Gamma(\wedge^{\bullet}A),\wedge,\lie{}{})$ such that $\delta^2=0$.
This is also called a \emph{differential Gerstenhaber algebra \cite{Xu:bv}}.
\end{exercise}

\begin{rmk}
Given a Lie bialgebroid $(A,\delta)$, there is a natural Lie algebroid structure on $A^*$ defined as follows.
\begin{enumerate}
\item The anchor map $\rho_*:A^*\to TM$ is given by \[ \ip{\rho_*\xi}{f}=\ip{\xi}{\delta f}, \quad \forall f \in \cty(M) .\]
\item The bracket $\lie{}{}$ is given by \begin{equation} \ip{\lie{\xi}{\eta}}{X}=(\delta X)(\xi,\eta)+(\rho_* \xi)\ip{X}{\eta}-(\rho_* \eta)\ip{X}{\xi}, \label{above} \end{equation}
for all $\xi,\eta\in\Gamma(A^*)$ and $X\in\Gamma(A)$.\end{enumerate}
Indeed, equivalently, a Lie bialgebroid is a pair of Lie algebroids $(A,A^*)$ such that \[ \delta\lie{X}{Y}=\lie{\delta X}{Y}+\lie{X}{\delta Y}, \quad \forall X,Y\in\Gamma(A) ,\] where $\delta:\Gamma(A)\to\Gamma(\wedge^2 A)$ is defined by the above equation \eqref{above}.
\end{rmk}

\begin{rmk} If $(A,A^*)$ is a Lie bialgebroid, then $(A^*,A)$ is also a Lie bialgebroid, called the dual Lie bialgebroid.
\end{rmk}

\begin{examples}
\begin{enumerate}
\item If $\pi$ is a \textbf{Poisson tensor} on $M$, then $A=TM$ with $\delta=\lie{\pi}{\cdot}:\mathfrak{X}^*(M)\to\mathfrak{X}^{*+1}(M)$ is a Lie bialgebroid.
In this case, $A^*=T^*M$ is the canonical cotangent Lie algebroid \cite{MackenzieX:1994}.
\item The \textbf{dual} to the previous one: $A=T^*M$, the cotangent Lie algebroid of a Poisson manifold $(M,\pi)$, together with $\delta^*=d_{\DR}:\Omega^*(M)\to \Omega^{*+1}(M)$ \cite{MackenzieX:1994}.
\item \textbf{Coboundary Lie bialgebroid.} Take $A$ a Lie algebroid admitting some $\Lambda\in\Gamma(\wedge^2A)$ satisfying \[ \mathcal{L}_X\lie{\Lambda}{\Lambda}=0, \forall X\in\Gamma(A) .\] Let $\delta=\lie{\Lambda}{\cdot}:\Gamma(\wedge^*A)\to \Gamma(\wedge^{*+1}A)$. Then $(A,\delta)$ defines a Lie bialgebroid \cite{LiuXu:1996}.
\item \textbf{Dynamical $r$-matrix.} 
\cite{EV, LiuXu:1999} 
Consider the Lie algebroid $A=T\mathfrak{h}^*\oplus\mathfrak{g}\to\eta$ where $\mathfrak{h}$ is an Abelian subalgebra of $\mathfrak{g}$ and the Lie algebroid structure on $A$ is the product Lie algebroid. Choose a map $r:\mathfrak{h}^*\to\wedge^2\mathfrak{g}$ and consider it as a element $\Lambda$ of $\Gamma(\wedge^2 A)$. Then
$\mathcal{L}_X\lie{\Lambda}{\Lambda}=0$ if, and only if, \[ \sum h_i\wedge\tfrac{dr}{d\lambda_i}+\tfrac{1}{2}\lie{r}{r}\in{(\wedge^3\mathfrak{g})}^{\mathfrak{g}} \] is a constant function over $\mathfrak{h}^*$. Here $\gendex{h_1, \dots,h_k}{}$ is a basis of $\mathfrak{h}$ and $(\lambda_1,\dots,\lambda_k)$ are the dual coordinates on $\mathfrak{h}^*$. \\
In particular, if $\mathfrak{g}$ is a simple Lie algebra and $\mathfrak{h}\subset \mathfrak{g}$ is a Cartan subalgebra, one can take \[ r(\lambda)=\sum_{\alpha\in\Delta^+}\frac{\lambda_{\alpha}}{(\alpha,\lambda)}e_{\alpha}\wedge f_{\alpha} \] or \[ r(\lambda)=\sum_{\alpha\in\Delta^+}\lambda_{\alpha} \coth{(\alpha,\lambda)}e_{\alpha}\wedge f_{\alpha} ,\] where $(e_{\alpha},f_{\alpha},h_i)$ is a Chevalley basis.
\end{enumerate}
\end{examples}


\section{Universal lifting theorem and quasi-Poisson groupoids}

The inverse procedure, i.e.  integration of Lie bialgebroids
to Poisson groupoids, is much more tricky than the 
case of groups.  This question was completely
 solved in \cite{MackenzieX:2000}.
In this section, we will  investigate
the integration problem from  a more general perspective point
of view. In particular, we show that  the integration
problem indeed follows from a general principle 
--- the ``universal lifting theorem'' --- 
in the theory of Lie groupoids,
which in turn implies other integration results.

\subsection{k-differentials}\label{4.1}
Let $A\to M$ be a Lie algebroid. Then $(\oplus \gm (\wedge^*A ), \wedge,
[\cdot, \cdot ])$ is a Gerstenhaber algebra.

\begin{defn} [\cite{ILX}]
A k-differential on a Lie algebroid $A$ is a degree $k-1$ derivation
of the Gerstenhaber algebra $(\oplus \gm (\wedge^*A ), \wedge,
[\cdot, \cdot ])$. I.e., a linear operator
\[ \delta : \gm (\wedge^* A)\to \gm (\wedge^{*+k-1}A) \]
satisfying
\begin{equation*}
\begin{array}{l}
\delta (P\wedge Q)=(\delta P)\wedge Q+(-1)^{p(k-1)}P\wedge \delta
Q, \\[7pt] \delta \lcf P, Q\rcf =\lcf \delta P, Q\rcf +(-1)^{(p-1)(k-1)}
\lcf P, \delta Q\rcf ,
\end{array}
\end{equation*}
for all $P\in \gm (\wedge ^pA)$ and $Q\in\gm (\wedge ^q A)$ .
\end{defn}

The space of all multi-differentials  $\cala = \oplus_{k \geq 0}  \cala_k$
becomes a graded Lie algebra under the graded commutator:
\[ [\delta_1, \delta_2]=\delta_1\smalcirc \delta_2
-(-1)^{(k-1)(l-1)}\delta_2\smalcirc \delta_1, \]
if $\delta_1\in \cala_k$ and $\delta_2\in \cala_l$.

The following result can be easily checked directly

\begin{lem}
A $k$ differential on a Lie algebroid $A$ is
equivalent to a pair of linear maps:
$C^\infty (M)\stackrel{\delta}{\to}\gm (\wedge^{k-1}A)$
and $\gm (A)\stackrel{\delta}{\to}\gm (\wedge^k A)$ satisfying
\begin{enumerate}
\item $\delta (fg)=g(\delta f)+f(\delta g)$,
for all $f, g\in  C^{\infty}(M)$;
\item $\delta (fX)=(\delta f)\wedge X+f\delta X$,
 for all $f \in C^{\infty}(M)$ and $X\in \gm (A)$;
\item $\delta \lcf X, Y\rcf =\lcf \delta X, Y\rcf +\lcf X, \delta Y\rcf$,
for all $X, Y\in \gm (A)$.
\end{enumerate}
\end{lem}

Below is a list of basic  examples.

\begin{example}
When $A$ is a Lie algebra $\mathfrak{g}$,
then a k-differential $\Leftrightarrow \ \delta: \mathfrak{g}\to \wedge^k \mathfrak{g}$
is a Lie algebra 1-cocycle w.r.t. the adjoint action.
\end{example}

\begin{example}
$0$-differential  $\Leftrightarrow$ $\phi\in \gm (A^*)$ such that
$d_A \phi =0$ is a Lie algebroid 1-cocycle.
\end{example}

\begin{example}
$1$-differential  $\Leftrightarrow $ infinitesimal of Lie algebroid
automorphisms \cite{MackenzieX:1998}.
\end{example}

\begin{example}
 $P \in \gm (\wedge^k A)$, then $\ad(P) =\lcf P , \cdot
\rcf$ is clearly a $k$-differential, which is called the \emph{coboundary} $k$-differential associated to $P$.
\end{example}

\begin{example}
A Lie bialgebroid $\Leftrightarrow$ a 2-differential of square
$0$ on a Lie algebroid $A$.
\end{example}

\subsection{Multiplicative $k$-vector fields on a Lie groupoid $\gm$}

Let $\gm\gpd M$ be a Lie groupoid and $\Pi  \in \mathfrak X ^k
(\gm)$.
Define $F_\Pi \in C^{\infty}(T^\ast \gm\times _\gm \stackrel{(k)}{\ldots }\times _\gm T^\ast \gm)$ by
\[ F_\Pi (\mu^1, \ldots , \mu^k)=\Pi (\mu^1, \ldots , \mu^k) \]

\begin{defn} [\cite{ILX}]
$\Pi \in \mathfrak X ^k (\gm)$ is multiplicative if and only if
$F_\Pi$ is a 1-cocycle w.r.t. the groupoid
$T^\ast \gm\times _\gm \stackrel{(k)}{\ldots }\times _\gm T^\ast
\gm \gpd A^\ast \times _M \stackrel{(k)}{\ldots }\times _M A^\ast $
\end{defn}

\begin{rmk}
$\Pi$ is multiplicative $\Leftrightarrow $ the graph of multiplication
$\Lambda \subset \gm \times \gm \times\gm$ is coisotropic w.r.t.
$\Pi\oplus \Pi \oplus (-1)^{k-1}\Pi$.
\end{rmk}

\begin{example}
If $P\in \gm (\wedge^k A)$, then $\Vec{P}-\ceV{P}$ is multiplicative.
\end{example}

By ${\mathfrak X}^k_{mult}(\gm)$ we denote the space of all multiplicative
k-vector fields on $\gm$. And ${\mathfrak X}_{mult}(\gm )=
\oplus_k {\mathfrak X}^k_{mult}(\gm )$.

\begin{prop}
The vector space ${\mathfrak X}_{mult}(\gm )$ is closed under the Schouten brackets
and therefore is a graded Lie algebra.
\end{prop}

The main result of this section is the following

\newtheorem*{ULT}{Universal Lifting Theorem}

\begin{ULT}[Iglesias-Laurent-Xu \cite{ILX}] 
If $\gm$ is $\alpha$-connected and
$\alpha$-simply connected,
then
\[ \mathfrak{X}_{mult}(\gm ) \cong \mathcal{A} \]
as graded Lie algebras.
\end{ULT}

\begin{proof}[Sketch of proof] 
Given a $\Pi \in {\mathfrak X}^k_{mult}(\gm )$.
Using the coisotropic condition, one can prove
that for any $f\in C^{\infty}(M)$ and $X\in \gm (A)$,
$[\beta^*f, \Pi ]$ and $[\ceV{X}, \Pi]$ are left invariant.
Define
$C^\infty (M)\stackrel{\delta_\Pi}{\to}\gm (\wedge^{k-1}A)$
and $\gm (A)\stackrel{\delta_\Pi}{\to}\gm (\wedge^k A)$
by
\[ \ceV{\delta_\Pi f}=[\beta^* f, \Pi ], \ \ \
\ceV{\delta_\Pi X}=[\ceV{X}, \Pi ] \]
One easily checks that
\begin{enumerate}
\item  $\delta_\Pi  (fg)=g(\delta_\Pi  f)+f(\delta_\Pi  g)$,
for all $f, g\in  C^{\infty}(M)$.
\item $\delta_\Pi  (fX)=(\delta_\Pi  f)\wedge X+f\delta_\Pi  X$,
 for all $f \in C^{\infty}(M)$ and $X\in \gm (A)$.
\item $\delta_\Pi  \lcf X, Y\rcf =\lcf \delta_\Pi  X, Y\rcf +\lcf X, \delta_\Pi  Y\rcf$,
for all $X, Y\in \gm (A)$.
\end{enumerate}
Thus $\delta_\Pi $ is a $k$-differential.
Moreover  the relation
\[ [\delta_\Pi, \delta_{\Pi'}]= \delta_{[\Pi, \Pi']} \]
 implies that $\Phi: \Pi\to \delta_\Pi$ is a Lie algebra homomorphism.

It is simple to check that $\Phi$ is injective. This is
because 
\be
\delta_\Pi =0 & \Leftrightarrow &  L_{\Vec{X}}\Pi =0 \mbox{ and } \Pi |_{M}=0\\
& \Leftrightarrow &\Pi=0.
\ee

The surjectivity needs more work. The main  idea is that
when $\gm$ is $\alpha$-connected and $\alpha$-simply connected
\[ \gm \cong P(A)/\sim \]
where $P(A) $ stands for the space of $A$-paths and 
$\sim $ is an equivalence relation on $A$-paths 
called homotopy, see \cite{CF} for more details. 
The quotient space $P(A)/\sim $ can be considered as the moduli space of flat connections
over the  interval $[0, 1]$ with the ``structure group'' being
the groupoid $\gm$.
Then
\begin{align*} 
& \text{k-differential on $A$ } \\ 
\Rightarrow \quad & \text{linear k-vector field on $A$} \\ 
\Rightarrow \quad & \text{k-vector field on the  space of paths $[0, 1]\to A$ } \\
\Rightarrow \quad & \text{descends to a well-defined k-vector field on $P(A)/\sim$}
\qedhere \end{align*}
\end{proof}

\begin{example}
In the coboundary case, the correspondence between
k-differentials and k-vector fields on $\gm$ is
given by
\[ \delta =[\Lambda , \cdot ], \ \Lambda \in \gm (\wedge^* A)
\qquad \leftrightarrow \qquad \Pi=\ceV{\Lambda}-\Vec{\Lambda} .\]
\end{example}

\subsection{Quasi-Poisson groupoids}

\begin{defn}[\cite{Roy}]
A quasi-Lie bialgebroid consists of a 2-differential $\delta$ on
a Lie algebroid $A$ such that $\delta^2 =[\phi, \cdot ]$ for
some $\phi \in \gm (\wedge^3 A)$ satisfying $\delta \phi =0$.
\end{defn}

\begin{example} \textbf{Quasi-Lie bialgebroids associated to twisted Poisson structures} \cite{SW, CX}
\label{twisted}

Recall that a twisted Poisson manifold \cite{SW} consists of
a triple 
 $(M,\pi,\omega)$, where
$\pi \in \mathfrak{X}^2 (M)$ and $\phi\in\Omega^3(M)$
satisfy the condition
 $\lie{\pi}{\pi}=(\wedge^3\pi^{\sharp})\phi$ and
 $d\phi=0$.

Let $T^*M$ be equipped with the following  Lie algebroid structure:
the bracket on $\Omega^1 (M)$ is given by
\[ \lie{\xi}{\eta}=\mathcal{L}_{\pi^{\sharp}\xi}\eta-\mathcal{L}_{\pi^{\sharp}\eta}\xi -d\pi(\xi,\eta)+\phi(\pi^{\sharp}\xi,
\pi^{\sharp}\eta,\cdot) , \ \ \forall \xi, \eta \in \Omega^1 (M) ;\]
 the anchor is $\rho=\pi^{\sharp}$.
One easily checks that this is indeed a Lie algebroid,
which is  denoted by $(T^*M)_{\pi,\phi}$.
Define
\[ \cinf{M}\xrightarrow{\delta}\Omega^1(M)\xrightarrow{\delta}\Omega^2(M) \] by 
$\delta(f)=df,\; \forall f\in\cinf{M}$ and $\delta(\eta)=d\eta-\pi^{\sharp}\eta\ii\phi,\;\forall\eta\in\Omega^1(M)$. 
Then $\delta$ is a 2-differential such that $\delta^2=\lie{\phi}{\cdot}$
Hence $\big((T^*M)_{\pi,\phi},\delta,\phi\big)$ is a quasi-Lie bialgebroid.
\end{example}

Assume that $\gm$ is an $\alpha$-connected and $\alpha$-simply connected
Lie groupoid with Lie algebroid $A$. Now let's see what the universal
lifting theorem implies. First, note that $\delta ^2=\frac{1}{2}[\delta,
\delta ]$. Hence we have
\[ \begin{cases}
\delta & \rightsquigarrow \quad \pi \text{ is a multiplicative bivector field on } \Gamma \\
\delta^2=\lie{\phi}{\cdot} & \rightsquigarrow \quad \tfrac{1}{2}\lie{\pi}{\pi}=\lvec{\phi}-\rvec{\phi} \\
\delta\phi=0 & \rightsquigarrow \quad [\pi,\lvec{\phi}]=0
\end{cases} \]

This motivates the following

\begin{defn}
A quasi-Poisson groupoid is a Lie groupoid $\gm\toto M$ with a
multiplicative bivector field $\Pi$  on $\gm$ and
$\phi \in \gm (\wedge^3 A)$  satisfying the identities  on the right
hand side above.
\end{defn}

The universal lifting theorem thus implies the following

\begin{thm}[\cite{ILX}]
If $\gm$ is $\alpha$-connected and
$\alpha$-simply connected Lie groupoid with Lie algebroid
$A$, then there is a
bijection between quasi-Lie bialgebroids $(A, \delta, \phi )$
and quasi-Poisson groupoids $(\Gamma\rightrightarrows M,\pi,\phi)$.
\end{thm}

Considering various special cases, we are led to a number
of  corollaries as below. In particular, we recover the
theorem of Mackenzie-Xu regarding the integration of
Lie bialgebroids.

When $\phi=0$, we have the following

\begin{cor}[\cite{MackenzieX:2000}]
If $\gm$ is $\alpha$-connected and
$\alpha$-simply connected Lie groupoid with Lie algebroid
$A$, then there is a
bijection between Lie bialgebroids $(A, \delta)$
and Poisson groupoids $(\Gamma\rightrightarrows M,\pi)$.
\end{cor}

Another special case is that $A$ is a Lie algebra.

\begin{cor}[\cite{Kosmann:1991}]
There is a bijection between quasi-Lie bialgebras
and connected and
simply connected quasi-Poisson groups.
\end{cor}

In particular,

\begin{cor}[\cite{Drinfeld83a}]
There is a bijection between Lie bialgebras
and connected and simply connected Poisson groups.
\end{cor}

Now let $(M,\pi)$ be a Poisson manifold. Then $(T^*M,d)$ is
a Lie bialgebroid, which gives rise to an $\alpha$-connected and
$\alpha$-simply connected Poisson groupoid
$(\gm \toto M, \Pi)$ (assuming that $\gm$ exists). In this case,
$M$ is called integrable;
 see \cite{CF2} for precise conditions for a Poisson manifold 
to be integrable.
 One shows that indeed $\Pi$ is non-degenerate in
this case. Thus one obtains a symplectic groupoid. This leads
to the following

\begin{cor} 
[\cite{Karasev, Weinstein:1987, CDW, CattaneoF}]
There exists a 1-1 correspondence between integrable  Poisson manifolds
 and $\alpha$-connected and $\alpha$-simply connected symplectic groupoids.

Every Poisson manifold can be realized as the base Poisson
manifold of  a (local) symplectic groupoid. In particular,
symplectic realizations always exist for any Poisson
manifold.
\end{cor}

Finally consider the quasi-Lie bialgebroid $((T^*M)_{\pi, \phi}, \delta, \phi)$ as in Example~\ref{twisted}.
Let $(\Gamma\rightrightarrows M,\Pi,\phi)$ be its corresponding 
$\alpha$-connected and $\alpha$-simply connected Poisson groupoid.
 One shows that $\Pi$ is non-degenerate. Let $\omega=\Pi^{-1}$.

Then one easily shows that 

\begin{enumerate}
\item $\omega$ is a multiplicative 2-form on $\Gamma$, i.e. the 2-form $(\omega,\omega,-\omega)$ vanishes when restricted to the graph of the multiplication
\item $d\omega=\alpha^*\phi-\beta^*\phi$
\end{enumerate}
That is, $(\Gamma\rightrightarrows M,\omega,\phi)$ is a twisted
symplectic groupoid.

Thus we obtain

\begin{cor}[\cite{CX}]
There exists a 1-1 correspondence between integrable twisted Poisson manifolds
 and $\alpha$-connected and $\alpha$-simply connected
twisted symplectic groupoids.
\end{cor}

\subsection{Symplectic  quasi-Nijenhuis groupoids}

Let $M$ be a smooth manifold, $\bivector$ a Poisson bivector field, 
and $\nijenhuis:TM\to TM$ a $(1,1)$-tensor.

\begin{defn}[\cite{MR1421686, MR1077465}]
The bivector field $\bivector$ and the tensor $N$ are said to be compatible if
\begin{gather*}
N\smalcirc\pi\diese=\pi\diese\smalcirc N^* \nonumber \\
\lie{\alpha}{\beta}_{\pi_N}=\lie{N^*\alpha}{\beta}_{\pi}
+\lie{\alpha}{N^*\beta}_{\pi}-N^*\lie{\alpha}{\beta}_{\pi}
\label{eq:newnir}
\end{gather*}
where $\pi_N$ is the bivector field on $M$ defined
by the relation $\pi_N\diese=\pi\diese\smalcirc N^*$ and
\begin{equation} \label{eq:pi}
\lie{\alpha}{\beta}_{\pi}:=\derlie{\pi\diese\alpha}(\beta)-
\derlie{\pi\diese\beta}(\alpha)-d\big(\pi(\alpha,\beta)\big),
\qquad \forall \alpha,\beta\in\df^1(M).
\end{equation}

\end{defn}

The $(1,1)$-tensor $N$ is said to be a Nijenhuis tensor if its
Nijenhuis torsion vanishes:
$$\lie{\nijenhuis X}{\nijenhuis Y}-\nijenhuis\big(\lie{\nijenhuis X}{Y}
+\lie{X}{\nijenhuis Y}-\nijenhuis\lie{X}{Y}\big)=0,\quad\forall X,Y\in\vf({M}).$$

In \cite{MR773513}, Magri and Morosi defined a 
 Poisson Nijenhuis manifold as  a triple $(M,\bivector,\nijenhuis)$,
where   $\pi$ is a Poisson bivector field, $N$ is a Nijenhuis tensor 
and $\pi$  and $N$ are compatible.

It is known that 
any Poisson Nijenhuis manifold $(M,\bivector,N)$ is endowed with a 
bi-Hamiltonian structure $(\pi, \pi_N )$, i.e. 
$$\lie{\bivector}{\bivector}=0,
\qquad\lie{\bivector}{\bivectorN}=0,
\qquad\lie{\bivectorN}{\bivectorN}=0.$$

Similarly, one can define Poisson quasi-Nijenhuis manifolds.

Let $i_N$ be the degree 0 derivation of $(\Omega^{\bullet}(M),\wedge)$ defined by
$$(i_N\alpha)(X_1,\cdots,X_p)=\sum_{i=1}^p\alpha(X_1,\cdots,NX_i,\cdots,X_p), \quad\forall\alpha\in\Omega^p(M).$$

\begin{defn} \label{def:42}
A Poisson quasi-Nijenhuis manifold is a quadruple $(M,\bivector,\nijenhuis,\ttf)$, 
where $\bivector\in\vf^2(M)$ is a Poisson bivector field, $N:TM\to TM$ is a $(1,1)$-tensor 
compatible with $\bivector$, and $\ttf$ is a closed 3-form on $M$ such that
$$\lie{\nijenhuis X}{\nijenhuis Y}-\nijenhuis\big(\lie{\nijenhuis X}{Y}+\lie{X}{\nijenhuis Y}
-\nijenhuis\lie{X}{Y}\big)=\bivector\diese(\interieur{X\wedge Y}\ttf), \quad\forall X,Y\in\vf(M)$$ 
and $i_N\ttf$ is closed.
\end{defn}

Defining a bracket 
$\lieN{\cdot}{\cdot}$ on $\XX(M)$ by $$[X,Y]_N =[NX,Y]+[X,NY]-N[X,Y],\quad\forall X,Y\in\XX(M)$$
as in \cite{MR1421686}, and considering $N:TM\to TM$ as an anchor map, we 
obtain a  new Lie algebroid structure on $TM$, denoted $(TM)_N$.
Its Lie algebroid cohomology differential $d_N: \Omega^{\bullet}(M)
\to \Omega^{\bullet+1}(M)$ is given by \cite{MR1421686}:
\begin{equation} d_N=[i_N, d]=i_N\smalcirc d-d \smalcirc i_N \label{eq:tgif} .\end{equation}

The following proposition extends a result of Kosmann-Schwarzbach \cite[Proposition 3.2]{MR1421686}.

\begin{prop}
\label{prop:qN}
The quadruple $(M,\bivector,\nijenhuis,\ttf)$ is a 
Poisson quasi-Nijenhuis manifold if, and only if,  
$\big((T^*M)_\pi, d_N,\ttf\big)$ is a quasi Lie bialgebroid 
and $\phi$ is a closed $3$-form.

In particular, the triple $(M,\bivector,\nijenhuis)$
 is a Poisson Nijenhuis manifold if, and only if,
$(T^*M)_\pi, d_N)$ is a  Lie bialgebroid.
\end{prop}

We now turn our attention to the particular case where the Poisson bivector field $\pi$ is 
non-degenerate. 

\begin{prop} \label{prop:3.10}
\begin{enumerate}
\item Let $(M,\pi,N,\ttf)$ be a Poisson quasi-Nijenhuis manifold. Then, 
\begin{equation} \lie{\pi}{\pi_N}=0, \label{eq:pqn2} \end{equation} and
\begin{equation} \lie{\pi_N}{\pi_N}=2 \pi\diese(\ttf) 
. \label{eq:pqn1} \end{equation}
\item Conversely, assume that $\pi\in\vf^2(M)$ is a non-degenerate Poisson bivector field, $N:TM\to TM$ is a $(1,1)$-tensor and $\ttf$ is a closed 3-form. If
 Eqs. \eqref{eq:pqn2}-\eqref{eq:pqn1} are satisfied, then $(M,\pi,N,\ttf)$ is a Poisson quasi-Nijenhuis manifold.
\end{enumerate}
\end{prop}

The following lemma is useful in characterizing
Poisson Nijenhuis structures in terms of Lie
bialgebroids.

\begin{lem}
  \label{lem:pnij}
  Let $(M,\pi)$ be a Poisson manifold.
  A Lie bialgebroid $((T^*M)_\pi,\delta)$
  is induced by a Poisson Nijenhuis structure if and only if
  $\lie{\delta}{d}=0$,
  where $d$ stands for the de Rham differential.
  \end{lem}

Given a Poisson Nijenhuis manifold
$(M,\pi,N)$, then $((T^*M)_\pi,d_N)$ is a Lie
bialgebroid.
 Assume that $(T^*M)_{\pi}$ is integrable,
and $(\gm\toto M,\tomega)$ is its $\alpha$-connected
and $\alpha$-simply connected symplectic groupoid.
Since $d_N^2=0$ and $\lie{d_N}{d}=0$, the universal lifting theorem implies that 
$d_N$ corresponds to a multiplicative Poisson bivector field $\tpi_\tN$ on $\gm$ such that 
$\lie{\tpi_\tN}{\tpi}=0$,
where $\tpi$ is the Poisson tensor on $\gm$ inverse to $\tomega$. 
Let $\tN=\tpi_\tN\diese\smalcirc\tomega\bemol:T\gm\to T\gm$.
Then it is clear that $\tN$ is a multiplicative
$(1,1)$-tensor, and the triple  $(\gm,\tomega,\tN)$
forms what is called    a {\em  symplectic Nijenhuis groupoid}.

\begin{defn}
A symplectic Nijenhuis groupoid is a symplectic groupoid
$(\gm \toto M, \tomega)$ equipped with a  multiplicative
$(1, 1)$-tensor $\tN: T\gm \to T\gm$ such that
$(\gm, \tomega, \tN )$  is a  symplectic Nijenhuis structure.
\end{defn}

Thus we have the following:

\begin{thm} [\cite{PqN}]
\label{thm:sym-nij}
\begin{enumerate}
\item The unit space of a symplectic Nijenhuis groupoid is a Poisson Nijenhuis manifold.
\item Every integrable Poisson Nijenhuis manifold is the unit space of a unique $\alpha$-connected,  $\alpha$-simply connected symplectic Nijenhuis groupoid.
\end{enumerate}
\end{thm}

Here, by an integrable Poisson Nijenhuis manifold, we mean the corresponding
 Poisson structure is integrable, i.e. it admits an associated symplectic
groupoid.

The  above theorem can be generalized
to the quasi-setting.

\begin{defn}
A symplectic quasi-Nijenhuis groupoid is a symplectic groupoid
$(\gm\toto M,\tomega)$ equipped with a multiplicative $(1,1)$-tensor 
$\tN:T\gm\to T\gm$ and a closed $3$-form $\ttf\in\Omega^3(M)$ 
such that $\big(\gm,\tomega,\tN,\beta^*\ttf-\alpha^*\ttf\big)$ is a symplectic 
quasi-Nijenhuis structure.
\end{defn}

The following result is a generalization of Theorem \ref{thm:sym-nij}.

\begin{thm}[\cite{PqN}]
 \label{thm:sym-nij1}
\begin{enumerate}
\item The unit space of a symplectic quasi-Nijenhuis groupoid
is a Poisson quasi-Nijenhuis manifold. 
\item Every integrable Poisson quasi-Nijenhuis manifold $(M,\pi,N,\ttf)$ is the unit space of a 
unique $\alpha$-connected and $\alpha$-simply connected symplectic quasi-Nijenhuis groupoid
$\big(\gm\toto M,\tomega,\tN,\beta^*\ttf-\alpha^*\ttf\big)$.
\end{enumerate}
\end{thm}

\subsection{Quasi-Poisson groupoid associated to Manin quasi-triple}

As an example, in this subsection, we discuss an important
class of quasi-Poisson groupoids which are associated to Manin quasi-triples.
For details, we refer the reader to \cite{ILX}.

Recall that a Manin pair $(\mathfrak{d},\mathfrak{g})$
\cite{Drinfeld}  consists of an even dimensional Lie algebra
 $\mathfrak{d}$ with an invariant, non-degenerate, symmetric bilinear form of signature $(n, n)$, and a Lagrangian subalgebra $\mathfrak{g}$. Given a Manin pair $(\mathfrak{d},\mathfrak{g})$, let $(D, G)$ be its corresponding
group pair (i.e., both $D$ and $G$ are connected and simply connected with Lie algebra $\mathfrak{d},\mathfrak{g}$ respectively).
The group $D$ and, in particular $G\subset D$, acts on $D/G$ by left multiplication. This is called the dressing action. One can form the corresponding transformation groupoid $G\times D/G\rightrightarrows D/G$ whose Lie algebroid is $\mathfrak{g}\times D/G\to D/G$.

Assume that $\mathfrak{h}\subset\mathfrak{d}$ is an isotropic complement to $\mathfrak{g}$: \[ \mathfrak{d}=\mathfrak{g}\oplus\mathfrak{h} .\] This yields a quasi-Lie algebra $(\mathfrak{g},F,\phi)$. Here, $F:\mathfrak{g}\to\wedge^2\mathfrak{g}$ and $\phi\in\wedge^3\mathfrak{g}$.

Let $\lambda:T^*(D/G)\to \mathfrak{g}\times D/G$ be the dual map of the dressing action \[ \mathfrak{g}^*\times D/G\simeq \mathfrak{h}\times D/G \xrightarrow{\text{dressing}}T(D/G) .\]

Consider $A=\mathfrak{g}\times D/G\to D/G$. Define \[ \delta:\cty(D/G)\to\Gamma(A)=\cty(D/G,\mathfrak{g}):f\mapsto \lambda(df) .\]
and
\[ \xymatrix{ 
\Gamma(A) \ar[r]^{\delta} \ar@2{-}[d] & \Gamma(\wedge^2 A) \ar@2{-}[d] \\ 
\xi\in\cinf{D/G,\mfg} \ar[r] & \cinf{D/G,\wedge^2\mfg}\ni F(\xi) } \] 
for $\xi$ a constant function.

Extend $\delta$ to all $\Gamma(\wedge^* A)$ using Leibniz rule.

\begin{prop}
$\delta$ is a 2-differential on $A$ and $\delta^2=\lie{\phi}{\cdot}$.
Hence $(A,\delta,\phi)$ is a quasi-Lie bialgebroid.
\end{prop}

\begin{thm}
[\cite{ILX}] 
Assume that $(\mathfrak{d},\mathfrak{g},\mathfrak{h})$ is a Manin
quasi-triple. Then $(G\times D/G\rightrightarrows D/G,\pi,\phi)$ is a quasi-Poisson groupoid where
\[
\Pi\Big((\theta_g,\theta_s),(\theta'_g,\theta'_s)\Big)=  \Pi_G(\theta_g,\theta'_g) -\Pi_{D/G}(\theta_s,\theta'_s)  +\ip{\theta'_s}{\widehat{(L_g^\ast \theta_g)}} -\ip{\theta _s}{\widehat{(L_g^\ast \theta'_g)}}
,\]
with 
\[ \Pi_{D/G}(df,dg)=-\sum \hat{\epsilon}^i(f) \hat{e}_i(g), \qquad\forall f,g\in \cty(D/G) .\] 
Here $\gendex{e_1,\dots,e_n}{}$ is a basis of $\mathfrak{g}$ and $\gendex{\epsilon_1,\dots,\epsilon_n}{}$ the dual basis of $\mathfrak{h}$.
\end{thm}

Let $\mathfrak g$ be a Lie algebra endowed with a non-degenerate
symmetric bilinear form $K$. On the direct sum $\mathfrak
d=\mathfrak g\oplus \mathfrak g$ one can construct a scalar
product $(\cdot | \cdot )$ by
\[
((u_1,u_2)|(v_1,v_2))=K(u_1,v_1)-K(u_2,v_2),
\]
for $(u_1,u_2),(v_1,v_2)\in \mathfrak d$. Then $(\mathfrak d,
\Delta (\mathfrak g),\half \Delta _- (\mathfrak g))$ is a Manin
quasi-triple with \[ \phi (u, v, w)=\frac{1}{4} K(u, [v, w]) .\]
(As usual, $\Delta (v) = (v,v) $ is the diagonal map, while
$  \Delta _- (v)= (v, -v)   $).
In this case, $D=G\times G$ and $D/G\cong G$. Through this identification,
$G$ acts  on $G$ by conjugation. Thus we have the 
following

\begin{cor}\label{GxG}
Assume that $\mathfrak g$ is a Lie algebra  endowed
with a non-degenerate symmetric bilinear form $K$ and $G$ is its
corresponding connected and simply connected Lie group. Then the
transformation groupoid $G\times G\gpd G$, where $G$ acts on $G$ by
conjugation, together with
the
multiplicative bivector  field $\Pi$ on $G\times G$:
\[
\Pi (g,s) = \half \sum _{i=1}^n \ceV{e^2_i}\wedge \Vec{e^2_i}-
\ceV{e^2_i}\wedge \ceV{e^1_i}-\Vec{(Ad_{g ^{-1}}e_i)^2}\wedge
\Vec{e^1_i},
\]
the bi-invariant 3-form $\Omega: =\frac{1}{4}K(\cdot ,[\cdot ,\cdot ]_\mathfrak g)\in \wedge^3 {\mathfrak g}^*\cong \Omega^3 (G)^G$ on $G$,
is  a quasi-Poisson groupoid. Here $\{ e_i\}$ is an orthonormal
basis of ${\mathfrak g}$ and the superscripts
refer to the respective $G$-component.
\end{cor}

\begin{example}
Another example is the case when 
$(\mathfrak{d},\mathfrak{g},\mathfrak{h})$
 is a Manin triple, i.e., $\phi=0$. Then one obtains a 
Poisson groupoid $G\times D/G\toto D/G$, which is
a symplectic groupoid integrating the Poisson
manifold $(D/G, \pi_{D/G})$. If moreover $G$ is complete,
$D/G\cong G^*$ as a Poisson manifold. Thus one obtains
the  symplectic groupoid of Lu-Weinstein 
$G\times G^*\toto G^*$ \cite{LuW:1989}.
\end{example}

\subsection{Hamiltonian $\gm$-spaces}

In this subsection, we show that the quasi-Poisson spaces
with $D/G$-valued momentum maps in the sense of Alekseev
and Kosmann-Schwarzbach correspond exactly to Hamiltonian $\gm$-spaces
of quasi-Poisson groupoids $\gm$.
 
Let $\Gamma \gpd M$ be a Lie groupoid. Recall that a {\em
$\gm$-space} is a smooth manifold $X$ with a map $J:X\to M$,
called the {\em momentum map},
and an action
\[ \Gamma \times _M X=\{ (g,x)\in \Gamma \times X \, |\, \bet (g)=J(x)\}\to X,
\ \ \
(g,x)\mapsto g\cdot x \] satisfying
\begin{enumerate}
\item $J(g\cdot x)=\alp (g)$, for $(g,x)\in \Gamma \times _M X$;
\item $(g h)\cdot x=g\cdot (h\cdot x)$, for $g,h\in \Gamma $
and $x\in X$ such that $\bet (g)=\alp (h)$ and $J(x)=\bet
(h)$;
\item $\epsilon (J(x))\cdot x =x$, for $x\in X$.
\end{enumerate}
Hamiltonian $\Gamma$-spaces for Poisson groupoids were studied in
\cite{LWX}. For quasi-Poisson groupoids, one can introduce
Hamiltonian $\gm$-spaces in a similar fashion.

\begin{defn}
\cite{ILX}
Let $(\Gamma \gpd M, \Pi , \Omega )$ be a quasi-Poisson groupoid.
A {\em Hamiltonian $\gm$-space}  is a $\gm$-space $X$ with momentum
map $J:X\to M$ and a bivector field $\Pi _X\in \mathfrak X ^2(X)$
such that:
\begin{enumerate}
\item the graph of the action $\{ (g,x,g\cdot x) \,|\, J(x)=\bet (g)
\}$ is a coisotropic submanifold of $(\Gamma \times X \times X,\Pi
\oplus \Pi _X\oplus -\Pi _X)$;
\item $\half[\Pi _X ,\Pi _X]=\hat{\Omega}$, where the hat denotes the map
$\Gamma (\wedge ^3A)\to \mathfrak X ^3(X)$, induced by the
infinitesimal action of the Lie algebroid on $X$: $\Gamma (A)\to
\mathfrak X (X)$, $Y\mapsto \hat{Y}$ .
\end{enumerate}
\end{defn}

\begin{prop}
Let $(\Gamma \gpd M, \Pi , \Omega )$ be a quasi-Poisson groupoid.
If $(X,\Pi _X)$ is a Hamiltonian $\Gamma$-space with momentum map
$J,$ then $J$ maps $\Pi _X$ to $\Pi _M$, where $\Pi _M$ is  the bivector
field on $M$ as in Corollary~\ref{corollary1}.
\end{prop}

Consider the quasi-Poisson groupoid $\gm: G\times D/G\toto D/G$
associated to a quasi Manin triple.  It is simple to check that
$(X,\Pi_X)$ is a Hamiltonian $\gm$-space iff there is a 
$G$-action on $X$ and a map $J:X\to D/G$ satisfying

\begin{enumerate}
\item $\Phi _\ast (\Pi _G\oplus \Pi _X) = \Pi _X$;
\item $\half [\Pi _X,\Pi _X]=\Omega _X$,
where $\Omega _X\in \mathfrak X ^3(X)$ is defined using the map
$\wedge ^3\mathfrak g \to \mathfrak X ^3(X)$ induced by the
infinitesimal action;
\item $\Pi _X^\sharp (J^\ast \theta _s)=-(\lambda (\theta _s))_X,\text{
for }\theta _s\in T^\ast _sS$, where $G$ acts on $D/G$ by dressing action.
\end{enumerate}

It deserves to be noted that the latter 
 is exactly a \emph{quasi-Poisson space} with $D/G$-momentum
map in the sense of Alekseev and Kosmann-Schwarzbach \cite{AK-S}.
In summary we have

\begin{thm}
[\cite{ILX}]
If $\gm$ is  the quasi-Poisson groupoid $\gm: G\times D/G\toto D/G$
associated to a quasi Manin triple, then
there is a bijection between 
Hamiltonian $\gm$-spaces and
quasi-Poisson spaces with $D/G$ -momentum map.
\end{thm}


\begin{thebibliography}{99}

\bibitem{AMM}
Alekseev, A., Malkin A., and Meinrenken, E.,
Lie group valued moment maps, {\em J. Diff. Geom.}
{\bf 48} (1998), 445-495.

\bibitem{AK-S} A. Alekseev and Y. Kosmann-Schwarzbach, Manin Pairs and
moment maps, {\em J. Differential Geom.} {\bf 56} (2000) 133--165.



\bibitem{CattaneoF}
  Cattaneo, A., and  Felder, G.,
     Poisson sigma models and symplectic groupoids,
 {\em Prog. Math.} {\bf 198} (2001) 61--93,

\bibitem{CattaneoF2}
 Cattaneo, A., and  Felder, G.,
     A path integral approach to the Kontsevich quantization,
{\em Comm. Math. Phys.} {\bf 212}  (2000), 591-611.

\bibitem{CX}
A.~S. Cattaneo and P. Xu, Integration of twisted Poisson
structures, {\em J. Geom. Phys.} {\bf 49} (2004) 187--196.

\bibitem{CP} V. Chari and A. Pressley, Quantum groups, Cambridge
  University Press.

\bibitem{CDW} A. Coste, P. Dazord and A. Weinstein,
Groupo{\"\i}des symplectiques, Publications du D{\'e}partement de
Math{\'e}matiques de l'Universit{\'e} de Lyon, {I}, {\bf 2/A}
(1987) 1--65.

\bibitem{CF} M. Crainic and R.-L. Fernandes, 
 Integrability of Lie brackets.  {\em Ann. of Math.} (2)  {\bf 157 }
 (2003),  575--620.

\bibitem{CF2} M. Crainic and R.-L. Fernandes, Integrability of Poisson
brackets, {\em J. Differential Geometry}. 
{\bf 66} (2004),  71--137.

\bibitem{Drinfeld83a}
Drinfel'd, V.G., Hamiltonian structures on Lie groups, Lie bialgebras,
and the geometric meaning of the classical Yang-Baxter equations, {\em
Soviet Math. Dokl.} {\bf 27}, (1983), 68-71.

\bibitem{Drinfeld83b}
Drinfel'd, V.G.,  On constant quasiclassical solutions of the Yang-Baxter
quantum equation,
{\em Soviet Math. Dokl.} {\bf 28}, (1983), 667-671.


\bibitem{Drinfeld} V. Drinfeld, Quasi-Hopf algebras, {\it
Leningrad Math. J.} \textbf{1} (1990) 1419--1457.

\bibitem{Dubrovin}
Dubrovin, B.,
 Geometry of $2$D topological field theories, in
{\em Integrable systems
and quantum groups}, Lecture Notes in Math. {\bf 1620}, Springer,
Berlin (1996), 120--348.

\bibitem{EV}  Etingof,  P.,  and Varchenko, A.,  Geometry and
classification of solutions of the classical dynamical Yang-Baxter
equation, {\it Commun. Math. Phys.}, {\bf 192} (1998), 77-120.

\bibitem{ILX}  D. Iglesias Ponte, C. Laurent-Gengoux and P. Xu,
 Universal lifting theorem and quasi-Poisson groupoids, math.DG/0507396.

\bibitem{Karasev} Karasev, M. V.,
   Analogues of the objects of Lie group theory for nonlinear
            Poisson brackets, {\em Math. USSR-Izv.} {\bf 28} (1987),
   497-527.


\bibitem{KM} M.V. Karas\"{e}v and , V. P. Maslov,
Asymptotic and geometric quantization. (Russian)
{\em Uspekhi Mat. Nauk} {\bf 39} (1984),  115--173.

\bibitem{Kosmann:1991}
Y.~Kosmann-Schwarzbach, Quasi-big\`ebres de Lie et groupes de Lie
quasi-Poisson, {\em C.R. Acad. Sci. Paris} {\bf 312} S\'erie I
(1991) 391--394.

\bibitem{MR1421686}
Y.~Kosmann-Schwarzbach, \emph{The {L}ie bialgebroid of a {P}oisson-{N}ijenhuis manifold},
  Lett. Math. Phys. \textbf{38} (1996), no.~4, 421--428.

\bibitem{MR1077465}
Y.~Kosmann-Schwarzbach and F.~Magri, \emph{Poisson-{N}ijenhuis structures},
  Ann. Inst. H. Poincar\'e Phys. Th\'eor. \textbf{53} (1990), no.~1, 35--81.


\bibitem{Laurent} C. Laurent-Gengoux,
From Lie groupoids to resolutions of singularity. Applications to 
symplectic resolutions (I), arXiv:math/0610288.



\bibitem{LX} C. Laurent-Gengoux and P. Xu,
 Quantization of pre-quasi-symplectic groupoids and their Hamiltonian
 spaces. 
{\em Progr. Math.}, {\bf 232}, (2005)  423--454. 

\bibitem{Soibel} S. Levendorskiui  and Y. Soibelman,  Algebras of
  functions on compact quantum groups, Schubert cells and quantum
  tori.  {\em Comm. Math. Phys.}  {\bf 139}  (1991),   141--170

\bibitem{Lie} S. Lie, Theorie der transformationsgruppen
(Zweiter Abschnitt, unter Mitwirkung von Prof. Dr. Friederich Engel).
Teubner, Leipzig, 1890.

\bibitem{LWX} Z.-J. Liu, A. Weinstein and P. Xu, Dirac structures and
Poisson homogeneous spaces, {\em Comm. Math. Phys.} {\bf 192}
(1998) 121--144.


\bibitem{LiuXu:1996}
Liu, Z.-J., and Xu, P.,
Exact Lie bialgebroids and Poisson groupoids,  {\em Geom. and Funtional Anal.}
 {\bf 6 } (1996) 138-145.

\bibitem{LiuXu:1999} Liu, Z.-J., and  Xu, P., Dirac structures
and dynamical $r$-matrices,
{\it Ann.  Inst. Fourier} {\bf 51} (2001), 831-859.



\bibitem{Lu:1990} J.-H. Lu, Multiplicative and affine Poisson structures on
Lie groups, {\sl Ph.D. Thesis}, University of California at
Berkeley, (1990).

\bibitem{LuW:1990} J.-H. Lu and A. Weinstein,   Poisson Lie groups,
dressing transformations, and Bruhat decompositions, {\em J. Diff.
Geom.} {\bf 31} (1990) 501--526.

\bibitem{LuW:1989}
Lu, J.-H. and Weinstein, A.,
Groupo{\"\i}des symplectiques doubles des groupes de
           {L}ie-{P}oisson,
{\em C. R. Acad. Sc. Paris}
{\bf 309} (1989), 951-954.

\bibitem{Mackenzie}
K.  Mackenzie,  General theory of Lie groupoids and Lie algebroids.
{\em  London Mathematical Society Lecture Note Series} {\bf  213}
 Cambridge University Press, Cambridge, 2005.

\bibitem{MackenzieX:1994}
Mackenzie, K.C.H., and Xu, P.,
Lie bialgebroids and Poisson groupoids,
{\em Duke Math. J.} {\bf 73} (1994), 415-452.

\bibitem{MackenzieX:2000}
Mackenzie, K.C.H.,  and Xu, P.,
Integration of Lie bialgebroids,  {\em Topology} {\bf 39} (2000), 445-467.





\bibitem{MackenzieX:1998} K.~C.~H. Mackenzie and P. Xu,
 Classical lifting processes and multiplicative vector fields,
 {\em Quarterly J. Math.} {\bf 49} (1998) 59--85.

\bibitem{MR773513}
F.~Magri and C.~Morosi, {On the reduction theory of the {N}ijenhuis
  operators and its applications to {G}elfand-{D}iki equations},
  Proceedings of the IUTAM-ISIMM symposium on modern developments in analytical
  mechanics, Vol. II (Torino, 1982), {\bf 117}, 1983, pp.~599--626.

\bibitem{PqN} Mathieu Sti{\'e}non and Ping Xu,
 Poisson quasi-{N}ijenhuis manifolds,  {\em Comm. Math. Phys.}
{\bf  270} (2007), 709--725,




\bibitem{Roy}
Roytenberg, D.,
 Courant algebroids, derived brackets and even symplectic
supermanifolds, Ph. D. thesis, Berkeley, 1999.


\bibitem{STS}
Semenov-Tian-Shansky, M.~A.,   Dressing transformations and Poisson Lie group
actions,
{\em Publ. RIMS, Kyoto University} {\bf 21} (1985), 1237-1260.

\bibitem{SW}
 Severa, P., and   Weinstein, A.,
Poisson geometry with a 3-form background,
{\em Prog. Theor. Phys. Suppl.} {\bf  144} (2001), 145-154.



\bibitem{Weinstein83} A. Weinstein,  The local structure of Poisson
  manifolds.  {\em J. Differential Geom.}  {\bf 18}  (1983),  523--557.

\bibitem{Weinstein:1987} A. Weinstein, Symplectic groupoids and Poisson
 manifolds,  {\em Bull. AMS} {\bf 16} (1987), 101-104.

\bibitem{Weinstein:cois}
Weinstein, A.,
 Coisotropic calculus and {P}oisson groupoids,
 {\em J. Math. Soc. Japan} {\bf 40} (1988), 705-727.



\bibitem{Weinstein:91} A. Weinstein, Symplectic groupoids, geometric
  quantization, and 
irrational rotation algebras.281--290, {\em Math. Sci. Res. Inst. Publ.},
  {\bf 20},  (1991) 281--290.

\bibitem{Xu:1995}
P. Xu,  On Poisson groupoids, {\em International J. of Math.}
{\bf 6 } (1995), 101-124.

\bibitem{Xu:bv}
P. Xu,  Gerstenhaber algebras and BV-algebras in  Poisson geometry,
{\em Commun. Math. Phys.} {\bf 200} (1999), 545-560.

\bibitem{Xu:Dirac}
Xu, P.,
 Dirac submanifolds and Poisson involutions,  {\it Ann. Scient. Ec.  Norm.
  Sup.}    {\bf 36} (2003), 403--430.



\bibitem{Xu:JDG04}
P. Xu, 
Momentum maps and Morita equivalence, {\it J. Diff. Geom.}  {\bf 67}
(2004), 289--333. 

\bibitem{Zak1}
S. Zakrzewski, Quantum and classical pseudogroups. Part 1. Union pseudogroups and their quantization, {\em Commun. Math. Phys.}{\bf 134} (1990), 347--370.

\bibitem{Zak2}
S. Zakrzewski, Quantum and classical pseudogroups. Part 2. Differential and symplectic pseudogroups, {\em  Commun. Math. Phys.}{\bf 134} (1990), 371--395.

\end{thebibliography}
\end{document}